\documentclass[12pt,a4paper]{article}
\usepackage[cp1251]{inputenc}
\usepackage[english]{babel}
\usepackage{graphicx}
\usepackage{indentfirst}
\usepackage{amsmath}
\usepackage{amsthm}
\usepackage{caption2}
\usepackage{wrapfig}
\usepackage[left=2.00cm, right=2.00cm, top=2.00cm, bottom=2.00cm]{geometry}
\usepackage{amsfonts}
\usepackage{graphics}

\newcommand{\bs}[1]{\boldsymbol{#1}}
\newcommand{\wt}[1]{\widetilde{#1}}

\newcommand{\const}{\mathrm{const}\,}
\newcommand{\pder}[2]{\frac{\partial #1}{\partial #2}}
\newcommand{\pdder}[2]{\frac{\partial^2 #1}{\partial #2^2}}

\newcommand{\Pc}{\mathcal{P}}

\newtheorem{Remark}{Remark}

\title{The motion of an unbalanced circular foil in the field of a point source}
\author{Elizaveta M. Artemova$^{1*}$, Evgeny V. Vetchanin$^{2+}$}
\date{}
\begin{document}
	\maketitle
	
	{\small \centering
		
		$^1$Ural Mathematical Center, Udmurt State University, Universitetskaya 1, Izhevsk, 426034 Russia

		$^2$Kalashnikov Izhevsk State Technical University, Studencheskaya 7, Izhevsk 426069, Russia
		
		$^*$liz-artemova2014@yandex.ru\quad $^+$eugene186@mail.ru
		
	}

	{ \small
		
		\textbf{Abstract.} 
		
		Describing the phenomena of the surrounding world is an interesting task that has long attracted the attention of scientists. However, even in seemingly simple phenomena, complex dynamics can be revealed. In particular, leaves on the surface of various bodies of water exhibit complex behavior.
		This paper addresses an idealized description of the mentioned phenomenon. Namely, the problem of the plane-parallel motion of an unbalanced circular disk moving in a stream of simple structure created by a point source is considered.
		Note that using point sources, it is possible to approximately simulate the work of skimmers used for cleaning swimming pools.
		Equations of joint motion of the unbalanced circular disk and the point source. It is shown that in the case of a fixed source of constant intensity the equations of motion of the disk are Hamiltonian. In addition, in the case of a balanced circular disk the equations of motion are integrable. A bifurcation analysis of the integrable case is carried out. Using a scattering map, it is shown that the equations of motion of the unbalanced disk are nonintegrable. The nonintegrability found here can explain the complex motion of leaves in surface streams of bodies of water.

	
	}
	
	\tableofcontents

	\section{Introduction}

The study of the motion of point singularities constitutes a fairly large part of theoretical
hydrodynamics. In this field, one of the classical problems going back to the works of
Helmholtz \cite{Helmholtz_1858} and Kirchhoff \cite{Kirchhoff_1874} is that of the
motion of rectilinear parallel vortex filaments in an ideal fluid, which is also called
the problem of the motion
of $N$ point vortices on a plane. It should be noted that
the equations of motion are integrable for $N \leq 3$, and in the case $N \geq 4$ this
system is nonintegrable and can exhibit chaotic behavior \cite{Ziglin_1980, Borisov_et_al_2006}.
	
	The conception of a point vortex was applied in various problem statements. For example,
this model was used to address the problem of the stability of
polygonal vortex configurations on a plane \cite{Kurakin_Yudovich_2002}, on a sphere
\cite{Borisov_Kilin_2000}, in spaces of constant curvature \cite{Borisov_et_al_2018}, and in a
circle \cite{Kurakin_2004}. A wide review of research results on vortex dynamics
is presented in the monograph \cite{Borisov_Mamaev_2005}. The problem of the perturbation of the motion of two point vortices by an acoustic wave is dealt with in \cite{Vetchanin_Kazakov_2016}. A similar 
problem with addition of background shear flow is treated in \cite{Vetchanin_Mamaev_2017}.
	
	Equations similar to the equations of motion of point vortices arise
in the description of the motion of point singularities in a Bose\,--\,Einstein condensate
\cite{Torres_et_al_2011, Koukouloyannis_et_al_2014, Ryabov_Sokolov_2019} and hetons
\cite{Sokolovskiy_et_al_2020a, Sokolovskiy_et_al_2020b}, oceanic vortex structures
existing in a stratified fluid. The description of the motion of vortex rings
also involves equations similar to the equations of motion of point vortices
\cite{Blackmore_Knio_2000a, Blackmore_Knio_2000b, Borisov_et_al_2013, Borisov_et_al_2014}.
	
	The point vortex model is also brought to bear to derive finite-dimensional equations
describing the joint motion of vortex structures and a rigid body. The equations of motion
of a balanced circular foil in the presence of point vortices were obtained for the first
time in \cite{Ramodanov_2001, Ramodanov_2002, Shashikanth_et_al_2002}. The motion
of an unbalanced circular foil in the presence of point vortices has recently been
examined in \cite{Mamaev_Bizyaev_2021}. We note that the displacement of the center of mass
of the circular foil leads to the loss of additional symmetry and
gives rise to obstructions to integrability. In particular, the problem of the motion of
a balanced circular foil and a point vortex is integrable, and in the case of a circular foil
with a displaced center of mass the system becomes nonintegrable and can exhibit
chaotic behavior.
	
	We note that the model of a point vortex can be used for a phenomenological description
of the shedding of vortices from a sharp edge of a rigid body performing unsteady motion \cite{Mason_2003, Michelin_Smith_2009, Fedonyuk_Tallapragada_2019}.
	
	In addition to point vortices, one considers (fairly rarely) the motion of other
point singularities in a fluid, such as sources, vortex sources, and dipoles
\cite{Fridman_Polubarinova_1928, Smith_2011, Bizyaev_et_al_2014, Bizyaev_et_al_2016}.
In our recent paper \cite{Artemova_Vetchanin_2020}, we proposed a finite-dimensional model
describing the motion of a balanced circular foil in the field of a fixed source.
It was shown that the system is integrable and admits one unstable periodic solution. It was
shown that this solution can be stabilized by changing the intensity of the source
via feedback.
	
	In this paper we develop the study presented in \cite{Artemova_Vetchanin_2020}, and
consider the motion of an unbalanced circular foil in the field of a point source.
We derive equations of joint motion of motion of a source whose intensity depends on time
and of an unbalanced circular foil. We show that, in the case of a fixed source of
constant intensity, the equations of motion of the circular foil are Hamiltonian.
We perform a bifurcation analysis of the dynamics for the case of a balanced foil.
In the system under consideration we have found no stable compact trajectories, and
therefore we carry out its numerical investigation using a scattering map instead of the
traditional Poincar\'{e} map. Using such a map, we show that, in the case of an unbalanced
foil, the equations of motion are nonintegrable.
	
	The model presented in this paper can be considered as an idealized description of the motion of leaves in streams on the surface of water bodies. Unconditionally, natural flows have a complex structure that is significantly different from the structure of the flow created by a point source.	Nevertheless, even in the case of the considered simple flow, the derived equations of motion turn out to be nonintegrable. The nonintegrability found here can be used as an explanation of the complex motion of leaves in streams observed in nature.
	
	
	\section{A mathematical model}
	
	\subsection{The case of a moving source of variable intensity}
	
	Consider the plane-parallel motion of an unbalanced circular foil in the presence
of a point source in an unbounded volume of an ideal incompressible fluid. We assume that
the fluid performs potential noncircular motion and is at rest at infinity.
	
	To construct a mathematical model, we introduce the following notation for the system parameters:
	\begin{itemize}
		\itemsep=-2pt
		\item $m_c$ is the mass of the foil,
		\item $I_c$ is the central moment of inertia of the foil,
		\item $R$ is the radius of the foil,
		\item $d$ is the distance between the geometric center of the foil and the center of mass,
		\item $q$ is the intensity of the source, which is generally a given
function of time.
	\end{itemize}
	Since we consider the plane-parallel motion, we will assume the parameters
$m_c$, $I_c$ and $q$ to be related to the unit of length of the foil. 
	
\begin{figure}[h!]
		\centering
		\includegraphics{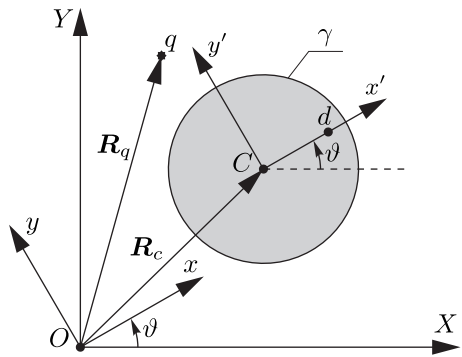}
		\caption{A schematic diagram of an unbalanced circular foil and a point source.}\label{fig.coord}
	\end{figure}

	To describe the motion of the system, we introduce three coordinate systems: a
fixed (inertial) system $OXY$, a moving system $Cx'y'$, rigidly attached to the foil,
and a coordinate system $Oxy$ rotating synchronously with the foil (see Fig. \ref{fig.coord}).
We will assume that the origin of the moving coordinate system $C$ is at the geometric center of
the foil and that the center of mass of the foil lies on the positive part of the axis $Cx'$.
	
	We will specify the position of the foil relative to the fixed coordinate system
by the radius vector of its geometric center $\bs R_c = (X_c,\, Y_c)$, and the orientation
of the foil, by the angle $\vartheta$ between the positive directions of the axes $OX$ and $Cx'$.
The position of the source relative to the fixed coordinate system will be specified
by the radius vector $\bs R_q = (X_q,\, Y_q)$. Thus, the configuration space of the system
is five-dimensional and
is $\mathcal{Q} = \bigl\{ \left( X_q,\, Y_q,\, X_c,\, Y_c,\, \vartheta \right)\ |\ (X_q - X_c)^2 + (Y_q - Y_c)^2 > R^2 \bigr\} \approx \mathbb{R}^3 \times \mathbb{T}^2$. The quantities defining the
configuration of the system are shown in Fig. \ref{fig.coord}.
	
	The motion of the system is determined by the interaction
of the foil and the source with the surrounding fluid. Since the motion of the fluid is
assumed to be potential, it can be completely described by a complex potential.
To construct the complex potential, to each point of the plane $OXY$ we associate the
complex number $Z = X+i Y$. Then the complex potential can be written as
	\begin{gather}\label{eq.W}
		W = - \frac{R^2\dot{Z}_c}{Z - Z_c} + \frac{q}{2 \pi} \ln(Z - Z_q) + \frac{q}{2 \pi} \ln \biggl( \frac{R^2}{Z - Z_c} - \overline{Z_q - Z_c} \biggr),
	\end{gather}
	where $Z_c = X_c + i Y_c$, $Z_q = X_q + i Y_q$ are complex-valued functions of time
which define the position of the center of the foil and the source, respectively.	
	\begin{Remark}
		The first and second terms in the expression \eqref{eq.W} are complex potentials
of the moving cylinder and the source, respectively \cite{Kochin}. The third term in the
expression \eqref{eq.W} arises by applying the Milne\,--\,Thomson theorem \cite{Miln_Thomson}
and describes the change of the flow created by the source, which occurs due to a circular
foil being added to it.
	\end{Remark}
	
	To construct the equations of motion, we introduce the following quantities:
	\begin{gather}\label{eq.ppp}
		\begin{gathered}
			P_x = m_c \dot{X}_c - m_c d \dot{\vartheta} \sin \vartheta, \quad P_y = m_c \dot{Y}_c + m_c d \dot{\vartheta} \cos \vartheta,\\
			P_{\vartheta} = m_c d ( - \dot{X}_c \sin \vartheta + \dot{Y}_c \cos \vartheta) + (I_c  + m_c d^2) \dot{\vartheta},
		\end{gathered}
	\end{gather}
	where $P_x$ and $P_y$ are the projections of the momentum of the foil onto the axes of
the fixed coordinate system and $P_\vartheta$ is the angular momentum of the foil relative
to its geometric center.

	The change in the momentum of the foil is defined by the principal vector of
pressure forces $\bs F = (F_x,\, F_y)$ acting on it from the fluid. The components
of the vector $\bs F$ can be calculated from Sedov's formula \cite{Sedov}:
	\begin{gather}\label{eq.SedovForce}
		F_x + i F_y = \overline{\frac{i \rho}{2} \oint\limits_{\gamma} \bigg( \frac{d W}{d Z} \bigg)^2 dZ} + \frac{d}{d t} \bigg( \rho \frac{d S Z_c}{d t} + i \rho \oint \limits_{\gamma} Z \frac{d W}{d Z} dZ \bigg).
	\end{gather}
	Here $\rho$ is the density of the fluid, $S = \pi R^2$ is the area of the foil, and
$\gamma$ is the boundary of the circular foil. Since the motion of the fluid is assumed
to be noncircular, the terms related to circulation are omitted in the formula
\eqref{eq.SedovForce}.
	
	Substituting the potential \eqref{eq.W} into the formula \eqref{eq.SedovForce}, we
obtain an explicit expression for the components of the principal pressure force vector:
	\begin{multline}\label{eq.forces}
		F_x + i F_y = \rho q R^2 \frac{(Z_q - Z_c)^2}{|Z_q - Z_c|^4} \dot{\overline{Z_q}} - \rho \dot{q} R^2 \frac{Z_q - Z_c}{|Z_q - Z_c|^2} + \\
		+ \frac{\rho q^2 R^2}{2 \pi} \frac{(Z_q - Z_c)}{(|Z_q - Z_c|^2 - R^2) |Z_q - Z_c|^2} - \rho \pi R^2 \ddot{Z}_c.
	\end{multline}
	The components $F_x$ and $F_y$ thus calculated are related to the unit of length of the cylinder.
The last term in the expression \eqref{eq.forces} coincides with the classical expression for
the force due to the effect of added masses \cite{Korotkin_2009} and acting on the
circular foil when it performs accelerated motion.
	
	\begin{Remark}
		We note that the calculation of the forces \eqref{eq.forces} can also be performed
in the real form
		\begin{gather*}
			F_x = -\oint\limits_{\gamma} p dy, \quad F_x = \oint\limits_{\gamma} p dx,
		\end{gather*}
		where the pressure $p$ is unambiguously calculated using the Cauchy\,--\,Lagrange
integral \cite{Kochin}. However, calculations in complex form with the aid of the formula
\eqref{eq.SedovForce} are simpler in practice.
	\end{Remark}

	Since the foil is circular, the pressure torque calculated relative to the geometric center
of the foil is zero. In this case, the angular momentum $P_\vartheta$ can change only
due to the rotation of the foil. Thus, the equations of motion of the foil take the form
	\begin{gather}\label{eq.general_eq}
		\dot{P}_x = F_x, \quad \dot{P}_y = F_y, \quad \dot{P}_{\vartheta} + d (P_x \cos \vartheta + P_y \sin \vartheta) \dot{\vartheta} = 0.
	\end{gather}
	
	We supplement equations \eqref{eq.general_eq} with the equations of motion of the source.
According to \cite{Fridman_Polubarinova_1928}, the velocity of the motion of the source
is equal to the velocity of the fluid at the point $Z = Z_q$ calculated
from the regular part of the potential~\eqref{eq.W} at the same point:
	\begin{gather}\label{eq.source}
		\dot{Z}_q = \overline{ \left(\frac{d W^*}{d Z}\right)} \bigg \vert_{Z = Z_q}, \quad W^* = W - \frac{q}{2 \pi} \ln(Z - Z_q).
	\end{gather}

	Equations \eqref{eq.ppp}, \eqref{eq.general_eq} and \eqref{eq.source} are a closed system and
completely describe the joint motion of the circular foil and the point source of
variable intensity in an ideal fluid.
	
	\subsection{The case of a moving source of constant intensity}
	
	In the case of a source of constant intensity ($\dot q = 0$) the equations of motion of
the foil \eqref{eq.general_eq} can be represented in the Lagrangian form
	\begin{gather}\label{eq.BodyLagr}
		\frac{d}{d t} \frac{\partial L}{\partial \dot{X}_c} - \frac{\partial L}{\partial X_c} = 0, \quad \frac{d}{d t} \frac{\partial L}{\partial \dot{Y}_c} - \frac{\partial L}{\partial Y_c} = 0, \quad \frac{d}{d t} \frac{\partial L}{\partial \dot{\vartheta}} - \frac{\partial L}{\partial \vartheta} = 0
	\end{gather}
	with the Lagrangian
	\begin{gather}\label{eq.L}
		L = T - U - \left( \bs A,\, \dot{\bs R}_c \right),
	\end{gather}
	where the following notation has been introduced:
	\begin{enumerate}
		\item $T$ is the kinetic energy of the foil and the fluid
		\begin{gather*}
			T = \frac{1}{2} \left( m_c + \rho \pi R^2 \right) \left( \dot{X}_c^2 + \dot{Y}_c^2 \right) + m_c d \left( -\sin\vartheta \dot{X}_c + \cos\vartheta \dot{Y}_c \right)\dot{\vartheta} + \frac{1}{2}\left( I_c + m_c d^2\right)\dot{\vartheta}^2,
		\end{gather*}
		\item $U$ is the scalar potential
		\begin{gather*}
			U = - \frac{\rho q^2}{4 \pi} \Bigl( \ln \bigl( (X_c - X_q)^2 + (Y_c - Y_q)^2 \bigr) - \ln\bigl( (X_c - X_q)^2 + (Y_c - Y_q)^2 - R^2\bigr) \Bigr),
		\end{gather*}
		\item $\bs A$ is the vector potential
		\begin{gather*}
			\bs A = -\frac{\rho q R^2}{(X_q - X_c)^2 + (Y_q - Y_c)^2} \left( \bs R_q - \bs R_c \right).
		\end{gather*}
	\end{enumerate}

	The right-hand sides of the equations of motion of the source \eqref{eq.source}
are also expressed in terms of the Lagrangian~\eqref{eq.L} as follows:
	\begin{gather}\label{eq.SrcLagr}
		\dot{X}_q = -\frac{1}{\rho q} \frac{\partial L}{\partial X_q}, \quad \dot{Y}_q = -\frac{1}{\rho q} \frac{\partial L}{\partial Y_q}.
	\end{gather}
	
	Thus, we have obtained a system of equations \eqref{eq.BodyLagr} and \eqref{eq.SrcLagr}
in a partially Lagrangian form which describes the joint motion of the
circular foil and the source of constant intensity.

	\subsection{The case of a fixed source of constant intensity}
	
	In the case of a fixed source of constant intensity ($\dot X_q = 0$, $\dot Y_q = 0$, $\dot q = 0$)
the equations of motion of the foil \eqref{eq.BodyLagr} can be represented in Hamiltonian form.
In this case, we will assume without loss of generality that the source
lies at the origin of the fixed coordinate system ($X_q = Y_q = 0$).
These conditions can always be enforced by making the change of variables
	\begin{gather*}
		X_c \to X_c - X_q,\quad Y_c \to Y_c - Y_q.
	\end{gather*}
	
	To reduce equations \eqref{eq.BodyLagr} to Hamiltonian form, we introduce new
generalized momenta $\Pi_x$, $\Pi_y$ and $\Pi_\vartheta$ as follows:
	\begin{gather}
		\Pi_x = \frac{\partial L}{\partial \dot{X}_c} = P_x + \rho \pi R^2 \dot{X}_c - A_x, \quad
		\Pi_y = \frac{\partial L}{\partial \dot{Y}_c} = P_y + \rho \pi R^2 \dot{Y}_c - A_y, \quad
		\Pi_{\vartheta} = \frac{\partial L}{\partial \dot{\vartheta}} = P_{\vartheta},\\
		A_x = \frac{\rho q R^2 X_c}{X_c^2 + Y_c^2},\quad A_y = \frac{\rho q R^2 Y_c}{X_c^2 + Y_c^2}\notag
	\end{gather}
	and perform the Legendre transformation
	\begin{gather}
		H = \Pi_x \dot{X}_c + \Pi_y \dot{Y}_c + \Pi_\vartheta \dot{\vartheta} - L \Big| _ {\dot{X}_c,\, \dot{Y}_c,\, \dot{\vartheta} \to \Pi_x,\, \Pi_y,\, \Pi_\vartheta} .
	\end{gather}
	Then the equations of motion of the circular foil \eqref{eq.BodyLagr} take the form
	\begin{gather}\label{eq.equation}
		\begin{gathered}
			\dot{X}_c = \frac{\partial H}{\partial \Pi_x}, \quad \dot{Y}_c = \frac{\partial H}{\partial \Pi_y}, \quad \dot{\vartheta} = \frac{\partial H}{\partial \Pi_{\vartheta}} \\
			\dot{\Pi}_x = -\frac{\partial H}{\partial X_c}, \quad \dot{\Pi}_y = -\frac{\partial H}{\partial Y_c}, \quad \dot{\Pi}_{\vartheta} = -\frac{\partial H}{\partial \vartheta},
		\end{gathered}
	\end{gather}
	where the Hamiltonian $H$ is given by the following expression:
	\begin{gather} \label{eq.Hamiltonian}
			H = \frac{1}{2} \big(\bs{\mathcal P}, \bold{Q}^{-1} \bs{\mathcal P} \big) - \frac{\rho q^2}{4 \pi} \Bigl( \ln \bigl( X_c^2 + Y_c^2 \bigr) - \ln\bigl( X_c^2 + Y_c^2 - R^2\bigr) \Bigr),\\
			\bs{\mathcal P} =
			\begin{pmatrix}
				\Pi_x + A_x \\
				\Pi_y + A_y\\
				\Pi_{\vartheta}
			\end{pmatrix}, \quad
			\bold{Q} =
			\begin{pmatrix}
				m_c + \rho \pi R^2 & 0 & -m_c d \sin \vartheta\\
				0 & m_c + \rho \pi R^2 & m_c d \cos \vartheta\\
				-m_c d \sin \vartheta & m_c d \cos \vartheta & I_c + m_c d^2
			\end{pmatrix}.\notag
	\end{gather}

	\begin{Remark}
		It has not been possible for us to represent in Hamiltonian form the complete system of
equations \eqref{eq.BodyLagr} and \eqref{eq.SrcLagr}, which describes the joint motion of
the circular foil and the point source of constant intensity. The form of the Lagrangian
\eqref{eq.L} and the form of equations \eqref{eq.SrcLagr} are obstructions
to such a representation.
		
		We note that, in a similar system governing the motion of the circular foil
in the presence of point vortices \cite{Bizyaev_Mamaev_2020},
representation in Hamiltonian form does turn out to be possible.
	\end{Remark}

	Equations \eqref{eq.equation} admit two first integrals: an energy integral coinciding
with the Hamiltonian \eqref{eq.Hamiltonian} and the integral of the angular momentum
	\begin{gather}\label{eq.Kint}
		K = \Pi_{\vartheta} + \Pi_y X_c - \Pi_x Y_c = \const,
	\end{gather}
	which is a consequence of the existence of a symmetry field
	\begin{gather}\label{eq.symField}
		\bs u = - Y_c \pder{}{X_c} + X_c \pder{}{Y_c} + \pder{}{\vartheta} - \Pi_y \pder{}{\Pi_x} + \Pi_x \pder{}{\Pi_y}.
	\end{gather}
	
	We next investigate the dynamics of the system and the properties of equations \eqref{eq.equation}
in the following two cases:
	\begin{itemize}
		\item[$-$] the case \textit{of a balanced} foil ($d = 0$);
		\item[$-$] the case \textit{of an unbalanced} foil ($d \neq 0$).
	\end{itemize}

	
	\section{A balanced cylinder}
	
	Consider the motion of a balanced ($d = 0$) circular foil in the field of a fixed source
of constant intensity. It turns out that in this case
a complete qualitative dynamics analysis can be carried out.
	
	If $d = 0$, the Hamiltonian \eqref{eq.Hamiltonian} takes the form
	\begin{gather}\label{eq.Hamiltonian_ball}
		H = H_{rot} + H_{trans},\\
		H_{rot} = \frac{1}{2} \frac{\Pc_\vartheta^2}{I_c},\label{eq.Hrot}\\
		H_{trans} = \frac{1}{2} \left( \frac{\Pc_x^2}{m} + \frac{\Pc_y^2}{m} \right) - \frac{\rho q^2}{4 \pi} \Bigl( \ln \bigl( X_c^2 + Y_c^2 \bigr) - \ln\bigl( X_c^2 + Y_c^2 - R^2\bigr) \Bigr),\label{eq.Htrans}
	\end{gather}
	where the following notation has been introduced:
	\begin{gather}
		\Pc_x = m \dot{X}_c, \quad \Pc_y = m \dot{Y}_c, \quad \Pc_\vartheta = m \dot{\vartheta}, \quad m = m_c + \rho \pi R^2.
	\end{gather}
	In this case, equations \eqref{eq.equation} decouple into two independent subsystems
governing the rotational motion
	\begin{gather}
		\begin{gathered}\label{eq.subsysRot}
			\quad \dot{\vartheta} = \pder{H}{\Pc_\vartheta} = \pder{H_{rot}}{\Pc_\vartheta},\quad \dot{\Pc}_\vartheta = - \pder{H}{\vartheta} = - \pder{H_{rot}}{\vartheta}
		\end{gathered}
	\end{gather}
	and the translational motion of the foil
	\begin{gather}
		\begin{gathered}\label{eq.subsysTrans}
			\dot{X}_c = \pder{H}{\Pc_x} = \pder{H_{trans}}{\Pc_x},\quad \dot{\Pc}_x = - \pder{H}{X_c} = - \pder{H_{trans}}{X_c} \\
			\dot{Y}_c = \pder{H}{\Pc_y} = \pder{H_{trans}}{\Pc_y},\quad \dot{\Pc}_y = - \pder{H}{Y_c} = - \pder{H_{trans}}{Y_c}
		\end{gathered}		
	\end{gather}
	with the Hamiltonians \eqref{eq.Hrot} and \eqref{eq.Htrans}, respectively.
	
	We first consider the subsystem \eqref{eq.subsysRot}. Since the
Hamiltonian \eqref{eq.Hrot} does not explicitly depend on $\vartheta$, this variable is
cyclic and the corresponding generalized momentum is a first integral:
	\begin{gather}\label{eq.Cint}
		C = \Pc_\vartheta = \const.
	\end{gather}
	On the fixed level set $C = c$ of the first integral  \eqref{eq.Cint}, the solution
of the subsystem \eqref{eq.subsysRot} can be represented as
	\begin{gather}
		\Pc_\vartheta = c = \const,\quad \vartheta (t) = \vartheta (0) + \frac{c}{I_c} t.
	\end{gather}

	Next, we consider the subsystem \eqref{eq.subsysTrans}, which can be investigated more conveniently by transforming to the polar coordinates:
	\begin{gather}\label{eq.polar1}
		X_c = s \cos \alpha,\quad Y_c = s \sin \alpha.
	\end{gather}
	In this case, we define the generalized momenta, which correspond to the coordinates $s$ and
$\alpha$, as follows:
	\begin{gather}\label{eq.polar2}
		\Pc _s = m \dot{s},\quad \Pc_\alpha = m s^2 \dot{\alpha}.
	\end{gather}
	The equations of motion \eqref{eq.subsysTrans} in the new variables
\eqref{eq.polar1} and \eqref{eq.polar2} remain canonical:
	\begin{gather}\label{eq.subsysTrans_polar}
		\dot{s} = \pder{H_{trans}}{\Pc_s},\quad \dot{\alpha} = \pder{H_{trans}}{\Pc_\alpha},\quad \dot{\Pc}_s = - \pder{H_{trans}}{s},\quad \dot{\Pc}_\alpha = - \pder{H_{trans}}{\alpha},
	\end{gather}
	where the Hamiltonian \eqref{eq.Htrans} takes the form
	\begin{gather}\label{eq.Htrans_polar}
		H_{trans} = \frac{1}{2 m} \left( \Pc_s^2 + \frac{\Pc_\alpha^2}{s^2} \right) + \frac{\rho q^2}{4 \pi} \ln \biggl( 1 - \frac{R^2}{s^2} \biggr).
	\end{gather}

	The Hamiltonian \eqref{eq.Htrans_polar} does not explicitly depend on the angle
$\alpha$. Hence, the generalized momentum $\Pc _\alpha$ is a first integral:
	\begin{gather}\label{eq.Fint}
		F = \Pc_\alpha = \const.
	\end{gather}

	\begin{Remark}
		In the case considered here ($d = 0$), using the expressions \eqref{eq.polar1}
and \eqref{eq.polar2}, the first integral \eqref{eq.Kint} can be represented as
		\begin{gather*}
			K = \Pc_\vartheta + \Pc_\alpha = C + F = \const.
		\end{gather*}
	\end{Remark}

	On the fixed level set $F = f$ of the first integral \eqref{eq.Fint} the equations for
$s$ and $\Pc_s$ decouple from the system \eqref{eq.subsysTrans_polar} and take the form
	\begin{gather}\label{eq.sPs}
		\dot{s} = \frac{\Pc_s}{m},\quad \dot{\Pc}_s = \frac{f^2}{m s^3} - \frac{\rho q^2 R^2}{2\pi s (s^2 - R^2)}.
	\end{gather}
	In this case, the evolution of the phase variable $\alpha$ is expressed by the quadrature
	\begin{gather*}
		\alpha(t) = \alpha(0) + \frac{f}{m} \int\limits_{0}^t \frac{d \tau}{s^2 (\tau)}.
	\end{gather*}
	This quadrature is necessary for reconstruction of the foil's motion relative to the fixed
coordinate system.

	Equations \eqref{eq.sPs} admit the first integral
	\begin{gather}\label{eq.Hred}
		H_{trans} = \frac{1}{2 m} \Pc_s^2 + U(s),\\
		U(s) = \frac{1}{2 m}\frac{f^2}{s^2} + \frac{\rho q^2}{4 \pi} \ln \biggl( 1 - \frac{R^2}{s^2} \biggr), \label{eq.potential}
	\end{gather}	
	which is a restriction of the Hamiltonian \eqref{eq.Htrans_polar}
to the fixed level set $F = f$ of the integral \eqref{eq.Fint}. Also,
equations \eqref{eq.sPs} have the involution
	\begin{gather}
		t \to -t,\quad s \to s,\quad \Pc_s \to -\Pc_s.
	\end{gather}

	\begin{Remark}\label{rem.3}
		The motions of the reduced system \eqref{eq.sPs} with $F = f > 0$ are identical with 
those with $F = -f$. Consequently, it suffices to analyze the dynamics for $f > 0$.
	\end{Remark}
	
	Equations \eqref{eq.sPs} are a Hamiltonian system with one degree of freedom.
In our analysis of this system we will use the classical approach of theoretical mechanics
which is based on the study of the level lines of the Hamiltonian (phase portrait) of the
system and bifurcations arising as the parameter values are varied \cite{Arnold, Bolsinov}.
	
	We first note that the Hamiltonian \eqref{eq.Hred} is defined for $s > R$ and,
regardless of the parameter values, possesses the following property:
	\begin{gather}
		\lim\limits_{s \to R + 0} H_{trans}(s, \Pc_s) = - \infty.\label{eq.singul}
	\end{gather}
	Also, the Hamiltonian \eqref{eq.Hred} is an even function of $\Pc_s$, i.e.,
$H_{trans}(s, -\Pc_s) = H_{trans}(s, \Pc_s)$, and hence the phase portrait of the system
will be symmetric relative to the line $\Pc_s = 0$.
It is seen from equations \eqref{eq.sPs} that the fixed points of the system
can lie only on the line $\Pc_s = 0$ and that they correspond to the critical points of the
potential \eqref{eq.potential}.
	
\begin{figure}[h]
	\centering
	\includegraphics{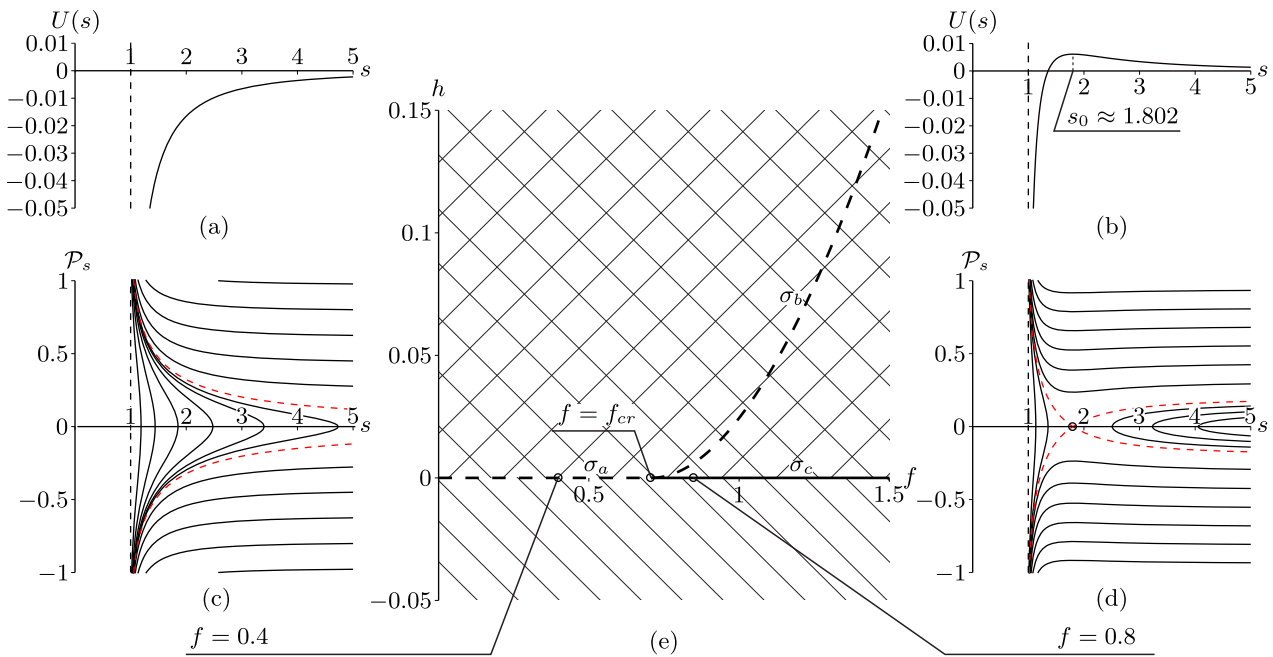}
	\caption{A typical view of a) the potential \eqref{eq.potential} for $f < f_{cr}$,
		b) the potential \eqref{eq.potential} for $f > f_{cr}$, c) the phase portrait of the system
		for $f < f_{cr}$, and d) the phase portrait of the system for $f > f_{cr}$.
		e) Bifurcation diagram. The parameter values: $m_c = 1$, $R = 1$, $q = 1$, $\rho = 1$.}
	\label{pic.diag}
\end{figure}

	Two qualitatively different cases can be singled out:
	\begin{itemize}
		\item[A.] The potential $U(s)$ (see Fig. \ref{pic.diag}b) has a maximum at the point
		\begin{gather}
			s_0 =  R |f| \sqrt{\frac{2 \pi}{2 \pi f^2 - \rho m q^2 R^2}}\label{eq.saddle}
		\end{gather}
		for $f > f_{cr}$, where
		\begin{gather*}
			f_{cr} = |q| R \sqrt{\frac{\rho m}{2 \pi}}.
		\end{gather*}
		The maximum point \eqref{eq.saddle} corresponds to the fixed saddle point of
the system \eqref{eq.sPs}:
		\begin{gather}\label{eq.saddle2}
			s = s_0,\quad \Pc_s = 0.
		\end{gather}
		A typical view of the phase portrait of the system
for $f > f_{cr}$ is shown in Fig. \ref{pic.diag}d. The red dashed lines in
Fig.~\ref{pic.diag}d indicate stable and unstable manifolds of the saddle point
\eqref{eq.saddle2}. The dashed line $s = R$ corresponds to the singularity \eqref{eq.singul}.
		
		\item[B.] If $f < f_{cr}$, the potential $U(s)$ is a monotonically increasing
function on the interval $s \in (R, +\infty)$ (see Fig. \ref{pic.diag}a). Hence, the system
has no fixed points. 		
		A typical view of the phase portrait of the system for $f < f_{cr}$ is shown in Fig.
\ref{pic.diag}c. The red dashed lines in Fig.~\ref{pic.diag}c denote the critical
trajectories separating different types of motion. The dashed line $s = R$
corresponds to the singularity \eqref{eq.singul}.
	\end{itemize}	
	
	It can be seen from the phase portraits (see Fig. \ref{pic.diag}c and \ref{pic.diag}d)
that all trajectories of the system \eqref{eq.sPs} either go to infinity or ``fall''
on a source, except for one trajectory corresponding to a fixed point when $f > f_{cr}$.
	
	A bifurcation diagram of the system on the plane of first integrals $(f,\, h)$, where
$h$ is a level set of the integral \eqref{eq.Hred}, is shown in Fig. \ref{pic.diag}e.
In view of Remark \ref{rem.3} we have presented only a part of the bifurcation diagram
which corresponds to the values $f > 0$. 	
	
\begin{figure}[h!]
		\centering
		\includegraphics{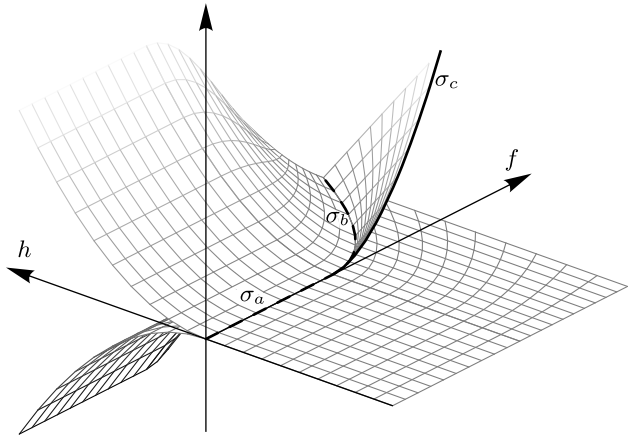}
		\caption{Bifurcation complex. The vertical axis in this picture has no physical meaning and is used merely
			for convenience of visualization of different leaves}\label{fig.leaves}
	\end{figure}

	If $f$ is fixed, to each value of $h > 0$ (the double hatched area in Fig. \ref{pic.diag}e)
there correspond two phase trajectories, and to the value of $h < 0$ (the single hatched area
in Fig. \ref{pic.diag}e), one phase trajectory. The bifurcation
diagram \ref{pic.diag}e corresponds to the bifurcation complex shown in Fig. \ref{fig.leaves}.
	
	The line $\sigma_a$ corresponds to critical trajectories of the system which arise
in the system for $f < f_{cr}$. The line $\sigma_b$ corresponds to the fixed point
\eqref{eq.saddle2} and to its stable and unstable manifolds. The line $\sigma_c$
corresponds to the boundary of the leaf.

	\begin{Remark}
		We note that the fixed point \eqref{eq.saddle2} of the reduced system \eqref{eq.sPs}
can be stabilized via feedback \cite{Artemova_Vetchanin_2020}. In this case, the equations
of motion cease to be Hamiltonian.
	\end{Remark}

	
	\section{An unbalanced cylinder}
	
	We now turn to an analysis of the motion of an unbalanced ($d > 0$) circular foil in
the field of a fixed source of constant intensity. It turns out that, as in the previous case
($d = 0$), examination of the system \eqref{eq.equation} in the case $d > 0$ by analytical 
methods reveals only a finite number of compact trajectories. Numerical experiments show that
the other trajectories of the system are noncompact, i.e., they either go to infinity or
``fall'' on the source. In addition, when $d > 0$, the equations of motion \eqref{eq.equation}
become nonintegrable.
	
	\subsection{The hypothesis of noncompactness of phase trajectories}\label{ssec.noncompact}
	
	To show the noncompactness of the phase trajectories of the system considered,
we single out the effective potential from the Hamiltonian \eqref{eq.Hamiltonian}.
To do so, we make the change of variables
	\begin{gather}\label{eq.polarTheta}
		x = X_c\cos\vartheta + Y_c\sin\vartheta,\quad y = -X_c\sin\vartheta + Y_c\cos\vartheta
	\end{gather}
	and express the Hamiltonian \eqref{eq.Hamiltonian} in terms of the derivatives of
the coordinates
$x$, $y$ and the angle $\vartheta$:
	\begin{gather}\label{eq.E0}
		\mathcal{H} = \frac{1}{2} \left( \bs{v},\, \wt{\bold{Q}} \bs{v} \right) + \frac{\rho q^2}{4\pi} \ln \left( 1 - \frac{R^2}{x^2 + y^2} \right),\\
		\wt{\bold{Q}} = \begin{pmatrix}
			m & 0 & 0 \\ 0 & m & m_c d \\ 0 &  m_c d  & I_c + m_c d^2
		\end{pmatrix},\quad \bs v = \begin{pmatrix} \dot{x} - y \dot{\vartheta} \\ \dot{y} + x \dot{\vartheta} \\ \dot{\vartheta} \end{pmatrix}.\notag
	\end{gather}
	The change of variables \eqref{eq.polarTheta} enables us to transform to the coordinate system
$Oxy$, which rotates synchronously with the foil, and to eliminate the angle variable
$\vartheta$ from the Hamiltonian. 
	
	We fix the level set $K = k$ of the first integral \eqref{eq.Kint} and express from it the derivative~$\dot{\vartheta}$:
	\begin{gather}\label{eq.dotTheta}
		\dot{\vartheta} = \frac{k + m y \dot{x} - (m x + m_c d) \dot{y}}{m(x^2 + y^2) + 2 m_c d x + I_c + m_c d^2}.
	\end{gather}
	Substituting the expression \eqref{eq.dotTheta} into the total energy \eqref{eq.E0},
we can write the following representation:
	\begin{gather}
		\mathcal{H} = \frac{1}{2} \left( \bs w,\, \bold{R}(x,\, y) \bs w \right) + U_e(x,\, y),\\
		U_e(x,\, y) = \frac{k^2}{2 \left( m(x^2 + y^2) + 2 m_c d x + I_c + m_c d^2 \right)} + \frac{\rho q^2}{4\pi} \ln \left( 1 - \frac{R^2}{x^2 + y^2} \right),\label{eq.potential2}
	\end{gather}
	where $\bs w = (\dot{x},\, \dot{y} )^T$, $\bold{R}(x,\, y)$ is
a symmetric positive definite matrix of fairly complex form, and $U_e(x,\, y)$
is the effective potential.
	
	Analysis of the expression \eqref{eq.potential2} shows that the critical points of the
effective potential can lie only in the plane $y = 0$. Thus, to find their coordinates,
it suffices to solve the equation
	\begin{gather}
		\pder{U_e(x,\, y)}{x} \bigg|_{y=0} = 0,
	\end{gather}
	which reduces to the fourth-degree equation:
	\begin{gather}
		a_4 x^4 + a_3 x^3 + a_2 x^2 + a_1 x + a_0 = 0,\label{eq.xEquation}\\
		\begin{gathered}
			a_4 = \frac{\rho q^2 R^2 m^2}{2\pi} - m k^2,\quad a_3 = \frac{2 \rho q^2 R^2 m m_c d}{\pi} - m_c d k^2,\\
			a_2 = \frac{\rho q^2 R^2}{2\pi} (2 m (I_c + m_c d^2) + 4 m_c^2 d^2) +  m R^2 k^2,\quad a_1 = \frac{2 \rho q^2 R^2 m_c d (I_c + m_c d^2)}{\pi} + m_c d R^2 k^2,\\
			a_0 = \frac{\rho q^2 R^2 (I_c + m_c d^2)^2}{2\pi}.
		\end{gathered}\notag
	\end{gather}
	
	\begin{Remark}
		In order to show that the critical points of the potential \eqref{eq.potential2} can
lie only in the plane $y = 0$, it suffices to make the change of variables
		\begin{gather*}
			x = r \cos\varphi,\quad y = r\sin\varphi
		\end{gather*}
		and to find the partial derivative with respect to the variable $\varphi$.
	\end{Remark}

	The analytical calculation of the roots of equation \eqref{eq.xEquation} leads to
cumbersome expressions whose analysis is difficult. Nonetheless,
numerical and analytical analyses of the function \eqref{eq.potential2} and
the coefficients of equation \eqref{eq.xEquation} reveal two level sets of the
integral \eqref{eq.Kint1}. When we pass through them, the surface of the
potential \eqref{eq.potential2} changes qualitatively:
	
	\begin{enumerate}
		\item The value $k_{inf}$ corresponds to the emergence of an inflection point on the curve
$U_e(x,\, 0)$. The coordinates of this point and the value $k_{inf}$ satisfy the system
of equations
		\begin{gather}
			\pder{U_e(x,\, y,\, k)}{x} \bigg|_{y=0} = 0,\quad \pdder{U_e(x,\, y,\, k)}{x} \bigg|_{y=0} = 0,
		\end{gather}
		which can be solved at least numerically.
		\item The value $k_{cr} = |q|R\sqrt{\dfrac{\rho m}{2\pi}}$ corresponds to the condition $a_4 = 0$.
In this case, equation \eqref{eq.xEquation} becomes cubic and necessarily has
one real root.
	\end{enumerate}

	
	Thus, one can single out five qualitatively different situations:
	\begin{enumerate}
		\item When $0 \leq |k| < k_{inf}$, the potential \eqref{eq.potential2}
has no critical points. A typical view of its surface and its profile in the plane $y = 0$
are shown
in Fig.\ref{fig.collage_a}.
		\begin{figure}[h!]
			\centering
			\includegraphics{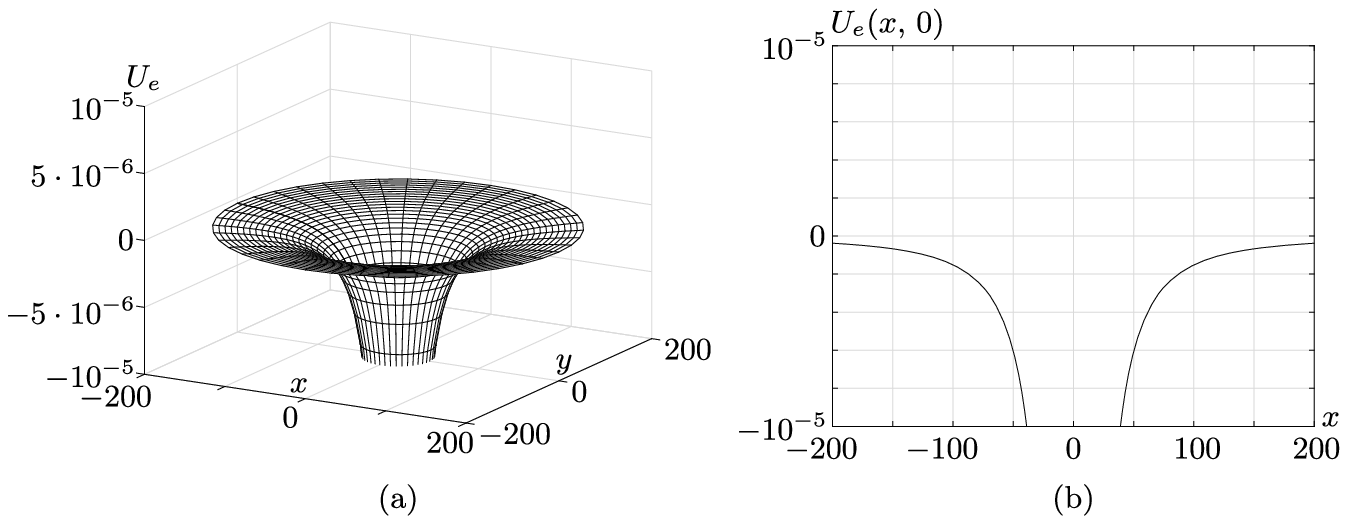}
			\caption{A typical view of (a) the surface of the potential and (b) the profile
of the function $U_e(x,\, 0)$ for $0 \leq |k| < k_{inf} \approx 0.81153$. The parameter
values: $m_c = 1$, $d = 0.1$, $R = 1$, $\rho = 1$, $q = 1$, $k \approx 0.73070$}\label{fig.collage_a}
		\end{figure}

		\item When $|k| = k_{inf}$, an inflection point arises on the curve $U_e(x,\, 0)$.
A typical view of the surface of the potential and its profile in the plane $y = 0$ are shown in Fig. \ref{fig.collage_b}.
		\begin{figure}[h!]
			\centering
			\includegraphics{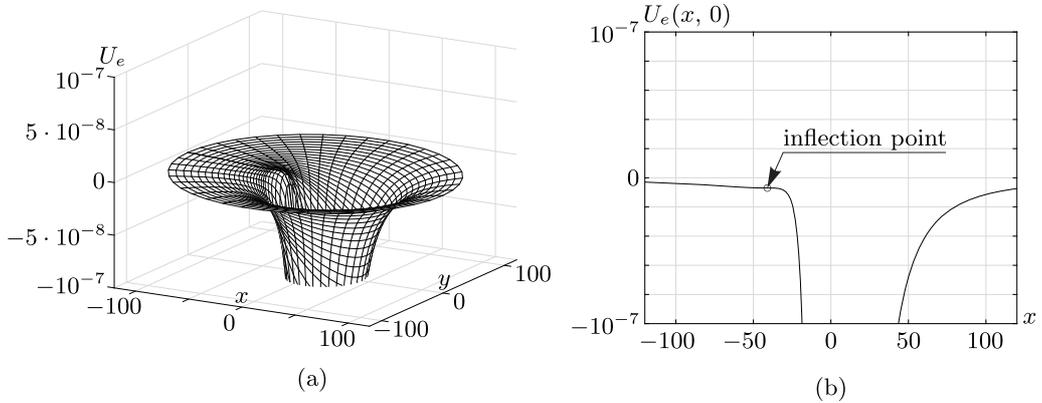}
			\caption{A typical view of (a) the surface of the potential and (b) the profile
of the function $U_e(x,\, 0)$ for $|k| = k_{inf} \approx 0.81153$. The parameter values: $m_c = 1$, $d = 0.1$, $R = 1$, $\rho = 1$, $q = 1$. The coordinate of the inflection point is $x_{inf} \approx -41.0591$}\label{fig.collage_b}
		\end{figure}
	

		\item If $k_{inf} < |k| < k_{cr}$, the potential \eqref{eq.potential2} has a
maximum point and a saddle point on the negative part of the axis $Ox$ (see Fig. \ref{fig.collage_c}).
Numerical analysis shows that the saddle point goes to minus infinity as $|k| \to k_{cr} - 0$.
		
		\begin{figure}[h!]
			\centering
			\includegraphics{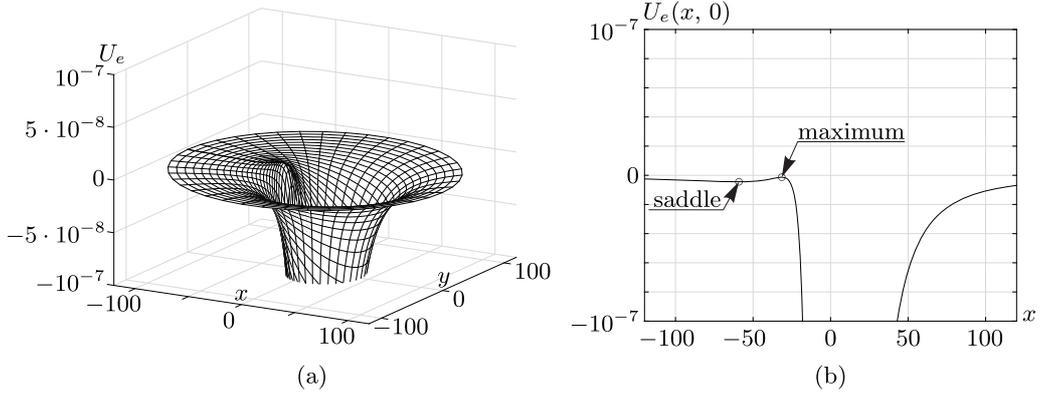}
			\caption{A typical view of (a) the surface of the potential and (b)
the profile of the function $U_e(x,\, 0)$ for $k_{inf} < |k| < k_{cr}$. The parameter values: $m_c = 1$, $d = 0.1$, $R = 1$, $\rho = 1$, $q = 1$, $k \approx 0.81156$, $k_{inf}\approx 0.81153$, $k_{cr} \approx 0.81188$. The coordinate of the saddle is $x_s \approx -59.09713$ and that of the
maximum point is $x_{max} \approx -31.45989$}\label{fig.collage_c}
		\end{figure}
	
		\item If $|k| = k_{cr}$, the potential \eqref{eq.potential2} has only a maximum point,
and the saddle point disappears. A typical view the surface of the potential and its profile
in the plane $y = 0$ are shown in
Fig. \ref{fig.collage_d}.
		\begin{figure}[h!]
			\centering
			\includegraphics{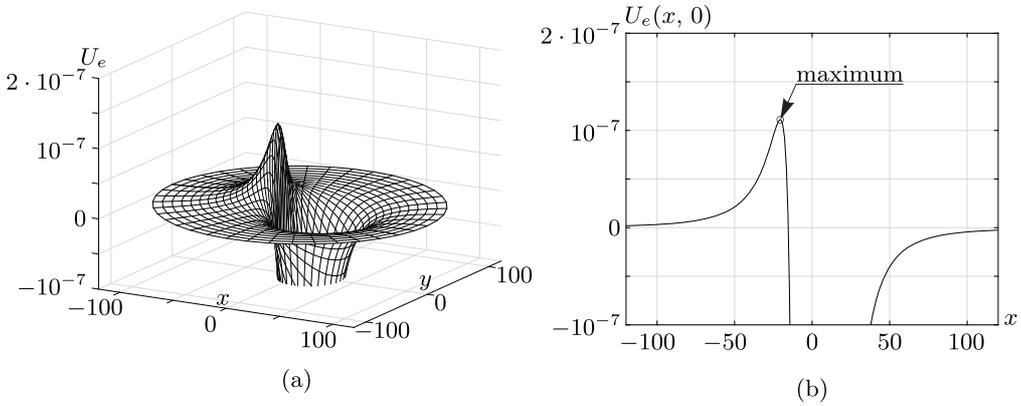}
			\caption{A typical view of (a) the surface of the potential and (b) the profile
of the function $U_e(x,\, 0)$ for $|k| = k_{cr} \approx 0.81188$. The parameter values: $m_c = 1$, $d = 0.1$, $R = 1$, $\rho = 1$, $q = 1$. The coordinate of the maximum point is $x_{max} \approx -20.54072$}\label{fig.collage_d}
		\end{figure}

		
		\item If $k_{cr} < |k|$, the potential \eqref{eq.potential2} has a maximum point
on the negative part of the axis $Ox$ and a saddle point on the positive part of the axis
(see Fig. \ref{fig.collage_e}). Numerical analysis
shows that the saddle point goes to plus infinity as $|k| \to k_{cr} + 0$.
		\begin{figure}[h!]
			\centering
			\includegraphics{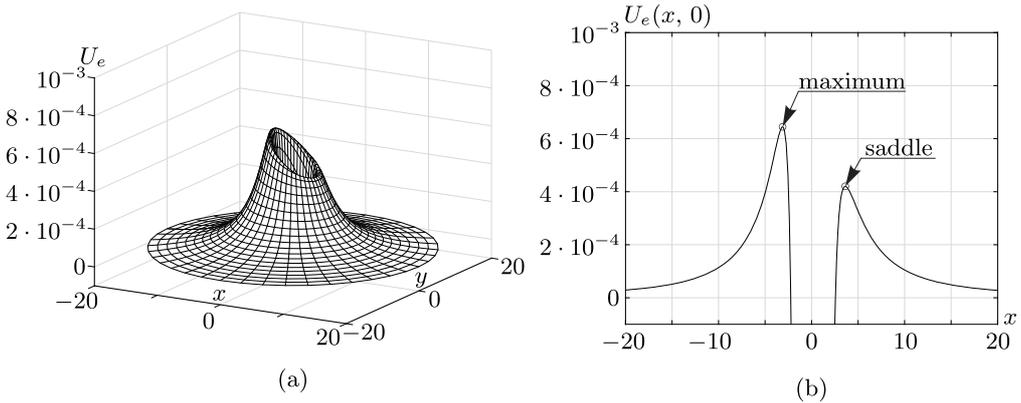}
			\caption{A typical view of (a) the surface of the potential and (b) the profile
of the function $U_e(x,\, 0)$ for $|k| > k_{cr} \approx 0.81188$. The parameter values: $m_c = 1$, $d = 0.1$, $R = 1$, $\rho = 1$, $q = 1$, $k \approx 0.86872$. The coordinate of the saddle
is $x_s \approx 3.62443$ and that of the maximum point is $x_{max} \approx -3.11904$}\label{fig.collage_e}
		\end{figure}
	
	\end{enumerate}
	
	The critical points of the effective potential $(x^*,\, 0)$ correspond to the
fixed points of the reduced equations of motion: $x = x^*$, $y = 0$, $\dot{x} = \dot{y} = 0$.
Here, by $x^*$ we mean the coordinate $x$ of any of the above-mentioned critical points.
In this case, the trajectory of the center of the foil is a circle:
	\begin{gather}\label{eq.circle}
		X_c(t) = x^* \cos \vartheta(t),\quad Y_c = x^* \sin\vartheta(t),\quad \vartheta(t) = \frac{k t}{m(x^*)^2 + 2 m_c d x^* + I_c + m_c d^2}.
	\end{gather}
	We note that, since the critical points of the potential are either maximum points or
saddle points, the motion in a circle \eqref{eq.circle} is unstable.
	
	We have seen that the effective potential has no minimum points. In addition,
numerical experiments show that the phase trajectories of the reduced system either go to
infinity or ``fall'' on the source. Thus, we can formulate the following hypothesis:
	\begin{enumerate}
		\item[ ] \textit{In the system under consideration there exist no compact
trajectories except for trajectories corresponding to the fixed points of the reduced system.}
	\end{enumerate}
	
\begin{figure}[h!]
		\centering
		\includegraphics[width=\linewidth]{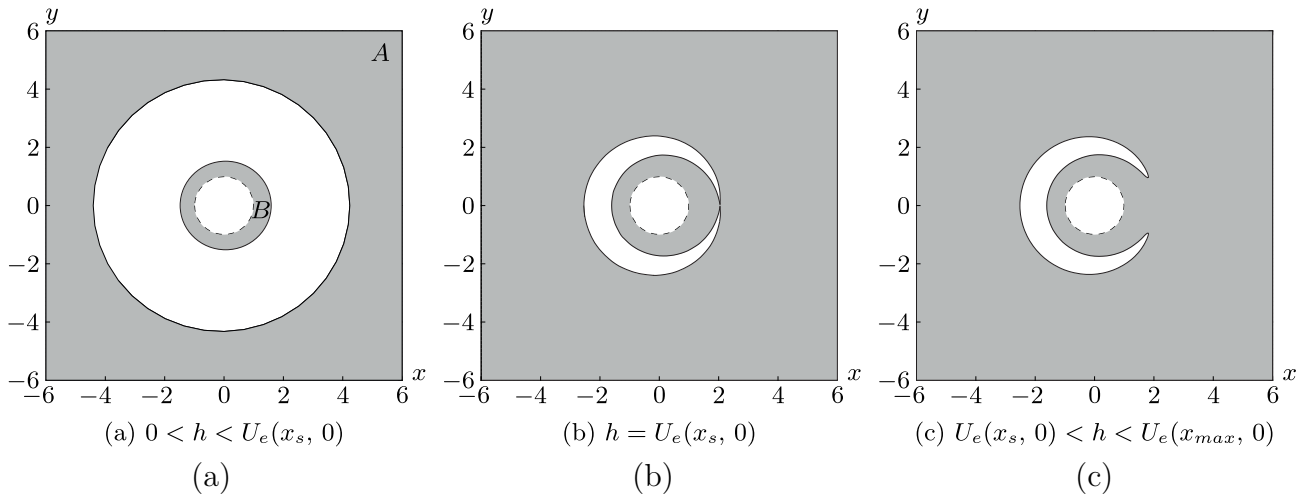}\\
		(a) \hspace{5cm} (b) \hspace{5cm} (c)
		\caption{Hill's regions for the case $|k| > k_{cr}$. The solid lines correspond to
the boundaries of Hill's regions, and the dashed line corresponds to a contact of the cylinder
with the source.}\label{fig.HillDomain}
	\end{figure}

	We now construct Hill's regions for the case $|k| > k_{cr}$. For $0 < h < U_e(x_s,\, 0)$
we have two Hill's regions divided by the potential barrier (see Fig. \ref{fig.HillDomain}a).
Thus, for the above-mentioned level sets of the energy integral the trajectories from region $A$
never fall on the source, and the trajectories from region $B$ never go to infinity.
When $h = U_e(x_s,\, 0)$, these regions come into contact with each other at the point $(x_s,\, 0)$, and
when $h > U_e(x_s,\, 0)$, there is only one Hill's region whose trajectories can
both go to infinity and fall on the source (see Fig. \ref{fig.HillDomain}b, c).
	
	\subsection{Nonintegrability of the equations of motion}
	
	We show that equations \eqref{eq.equation} are nonintegrable. One of the
classical numerical methods of investigating the nonintegrability of equations of
motion is the method of constructing a Poincar\'{e} map \cite{Kuznetsov_2006}.
However, such a map requires that the phase trajectory return to the secant. The analysis
made in subsection \ref{ssec.noncompact} has shown that almost all trajectories of the
system \eqref{eq.equation} are noncompact, except for a finite number of unstable periodic
trajectories, which makes the traditional Poincar\'{e} map inapplicable. Nonetheless,
the nonintegrability of equations \eqref{eq.equation} can be investigated
using a scattering map, which is an analog of the
Poincar\'{e} map for
systems with noncompact trajectories~\cite{Jung_1986}.
	
	We will consider the motion on the level sets of the integral \eqref{eq.Kint} $|k| > k_{cr}$ 
and the levels of energy corresponding to the case shown in Fig. \eqref{fig.HillDomain}a. 
The restrictions made are sufficient to ensure that the foil's trajectories 
coming from infinitely remote points will not ``fall'' on the source, but go back to infinity.

	To construct a scattering map, it is convenient to write the equations of motion 
\eqref{eq.equation} in the momenta relative to the rotating coordinate system $Oxy$. 
To do so, we make the change of
variables \eqref{eq.polarTheta}
	\begin{gather*}
		x = X_c\cos\vartheta + Y_c\sin\vartheta,\quad y = -X_c\sin\vartheta + Y_c\cos\vartheta
	\end{gather*}
	and denote the new momenta:
	\begin{gather}\label{eq.momentums}
		p_x = (\Pi_x+A_x)\cos\vartheta + (\Pi_y+A_y)\sin\vartheta,\quad p_y = -(\Pi_x+A_x)\sin\vartheta + (\Pi_y+A_y)\cos\vartheta,\quad p_\vartheta = \Pi_\vartheta.
	\end{gather}

	The change of variables \eqref{eq.polarTheta} and \eqref{eq.momentums} satisfies
the relations
	\begin{gather}\label{eq.varProp}
		\bs u(x) = \bs u(y) = \bs u(p_x) = \bs u(p_y) = \bs u(p_\vartheta) = 0,\quad  \bs u(\vartheta) = 1,
	\end{gather}
	where $\bs u$ is the symmetry field \eqref{eq.symField}. Due to \eqref{eq.varProp},
the variable $\vartheta$ becomes cyclic and the equations for $x$, $y$, $p_x$, $p_y$ and 
$p_\vartheta$ decouple from the complete
system of equations. In this case, the variable $\vartheta$ is necessary only 
to reconstruct the motion in the full phase space and need not be considered in the further 
analysis of integrability, see, e.g., \cite{Borisov_Mamaev_2015}.
	
	In the new variables the Hamiltonian \eqref{eq.Hamiltonian} becomes
	\begin{gather}
		H = \frac{1}{2} \left( \bs p,\, \wt{\bold Q}^{-1} \bs p \right) - \frac{\rho q^2}{4 \pi} \Bigl( \ln \bigl( x^2 + y^2 \bigr) - \ln\bigl( x^2 + y^2 - R^2\bigr) \Bigr),\\
		\bs p = \begin{pmatrix}
			p_x\\ p_y\\ p_{\vartheta}
		\end{pmatrix},\quad \wt{\bold{Q}} = \begin{pmatrix}
			m & 0 & 0 \\ 0 & m & m_c d \\ 0 & m_c d & I_c + m_c d^2
		\end{pmatrix}. \notag
	\end{gather}
	We recall that $m = m_c + \rho \pi R^2$.

	The nonzero Poisson brackets in the new variables have the form
	\begin{gather}
		\begin{gathered}
			\{ p_x,\, p_\vartheta \} = p_y,\quad \{ p_y,\, p_\vartheta \} = -p_x,\quad \{ x,\, p_x \} = \{ y,\, p_y \} = \{ \vartheta,\, p_\vartheta \} = 1,\\
			\{ x,\, p_\vartheta \} = y,\quad \{ y,\, p_\vartheta \} = -x.
		\end{gathered}
	\end{gather}
	
	The integral \eqref{eq.Kint} in the new variables preserves its form:
	\begin{gather}\label{eq.Kint1}
		K = p_{\vartheta} + p_y x - p_x y = \const.
	\end{gather}
	On the fixed level set $K = k$ of the first integral \eqref{eq.Kint1} the variable
$p_\vartheta$ can be eliminated from the equations of motion. In addition,
to analyze the dynamics, it is more convenient to make another change of variables:
	\begin{gather}
		x = r \cos\varphi,\quad y = r \sin\varphi, \quad p_x = p \cos\alpha,\quad p_y = p \sin\alpha.
	\end{gather}
	
	The reduced system of equations of motion can be represented as
	\begin{gather}\label{eq.subsys}
		\begin{gathered}
			\dot{r} = m^{-1} p \cos (\alpha - \varphi) - m^{-1}m_c d \sin\varphi \Omega,\quad
			\dot{\varphi} = \frac{p\sin(\alpha - \varphi)}{m r} - \frac{m_c d \cos\varphi}{m r} \Omega - \Omega,\\
			\dot{p} = - \frac{\rho q^2 R^2}{2\pi} \frac{\cos(\alpha - \varphi)}{r (r^2 - R^2)},\quad
			\dot{\alpha} = \frac{\rho q^2 R^2}{2 \pi} \frac{\sin(\alpha - \varphi)}{r (r^2 - R^2) p} - \Omega,
		\end{gathered}\\
		\Omega = \frac{k - r p \sin(\alpha - \varphi) - m^{-1} m_c d p \sin\alpha}{I_c + m_c d^2 - m^{-1} m_c^2d^2},\label{eq.dotTheta1}
	\end{gather}
	with the Hamiltonian
	\begin{gather}\label{eq.Ham}
		H = \frac{1}{2} \left( m^{-1} p^2 + (I_c + m_c d^2- m^{-1} m_c^2 d^2) \Omega^2 \right) + \frac{\rho q^2}{4\pi} \ln\left( 1 - \frac{R^2}{r^2}\right)
	\end{gather}
	and the Poisson bracket
	\begin{gather}
		\begin{gathered}
			\{ r,\, p \} = \cos(\varphi - \alpha),\quad \{ r,\, \alpha \} = \frac{\sin(\varphi - \alpha)}{p},\\
			\{ \varphi,\, p \} = -\frac{\sin(\varphi - \alpha)}{r},\quad \{ \varphi,\, \alpha \} = \frac{\cos(\varphi - \alpha)}{r p},\quad 
		\{\vartheta,\, K \} = 1.	
		\end{gathered}
	\end{gather}

	We see that the equations in the variables $r$, $\varphi$, $p$ and $\alpha$ decouple 
from the complete system and possess the first integral \eqref{eq.Ham}. Equation \eqref{eq.dotTheta1} 
and the first integral \eqref{eq.Kint1} are necessary to reconstruct the motion in the 
complete phase space of the system.
	
	We note that the terms with $\dot{\vartheta}$ arise in equations \eqref{eq.subsys} 
due to transformation to the moving coordinate system and lead to a drift in the system. 
In fact, we have the problem of a material point moving on the surface of a potential 
rotating with angular velocity $-\dot{\vartheta}$.
	
	From the form of equations \eqref{eq.subsys} and the Hamiltonian \eqref{eq.Ham} 
it can be seen that, as $r\to\infty$, the following holds:
	\begin{gather}\label{eq.asymptotic}
		\begin{gathered}
			b = r \sin(\alpha - \varphi) + m^{-1} m_c d \sin\alpha \to \const,\quad \dot{\vartheta} \to \const,\\
			\alpha - \varphi \to \pi n,\quad \alpha + \vartheta \to \const,\quad \varphi + \vartheta \to \const, \quad p \to \const.
		\end{gathered}
	\end{gather}
	Thus, as the foil recedes to a sufficiently great distance from the source, the foil 
rotates almost uniformly with angular velocity $\dot{\vartheta}$. In \eqref{eq.asymptotic} 
the quantity $b$ has the meaning of an impact parameter.

	In order to demonstrate the nonintegrability of equations \eqref{eq.subsys}, 
we construct a scattering map, as is done in \cite{Jung_1986}. This map is a composition 
of two maps: a direct map $\mathcal{F}$ and a feedback map $\mathcal{B}$.

\begin{figure}[h!]
		\centering
		\includegraphics{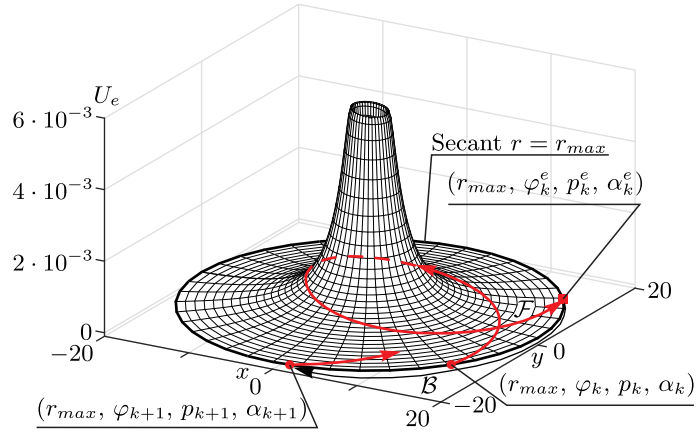}
		\caption{A schematic illustration of a scattering map. 
$\mathcal{F}$ is the direct map and $\mathcal{B}$ is the feedback map}\label{fig.scatteringMapScheme}
	\end{figure}

	\begin{enumerate}
		\item The direct map $\mathcal{F}$ consists in integrating equations 
\eqref{eq.subsys} with the initial condition $(r_{max},\, \varphi_k,\, p_k,\, \alpha_k)$. 
In this case, the hypersurface $r = r_{max}$ is a secant of the scattering map, and 
the value of momentum $p_k$ is calculated from a given level of energy and the values 
$r_{max}$, $\varphi_k$ and $\alpha_k$. We note that, since the Hamiltonian \eqref{eq.Ham} 
is a quadratic function of momentum $p$, there exist on each fixed level of energy 
two different values of momentum $p$. As $p_k$ we need to choose solely either
the largest or the smallest value of momentum.
		
		Integration of equations \eqref{eq.subsys} is performed until the phase trajectory
reaches again the secant $r = r_{max}$ (see Fig.\ref{fig.scatteringMapScheme}). 
We denote the end point as $(r_{max},\, \varphi_k^e,\, p_k^e,\, \alpha_k^e)$. 
It should be kept in mind that, regardless of the choice of the initial momentum, the 
value $p_k^e$ can correspond both to the largest and to the smallest
of the possible values of momentum.
		
		\item The feedback map (see Fig. \ref{fig.scatteringMapScheme}) 
consists in transforming the final values of the angle variables $\alpha_k^e$ and 
$\varphi_k^e$ to the new initial values as follows:
		\begin{gather}
			\alpha_{k+1} = \alpha_k^e,\quad \varphi_{k+1} = 2\alpha_k^e - \varphi_k^e + \pi.
		\end{gather}
	\end{enumerate}
	Applying alternately the direct map $\mathcal{F}$ and the feedback map $\mathcal{B}$, 
we obtain the following two-dimensional map:
	\begin{gather}\label{eq.scatteringMap}
		\mathcal{B}\mathcal{F} : \ (\varphi_k,\, \alpha_k) \to (\varphi_{k+1},\, \alpha_{k+1}).
	\end{gather}

\begin{figure}[h!]
	\centering
	\includegraphics{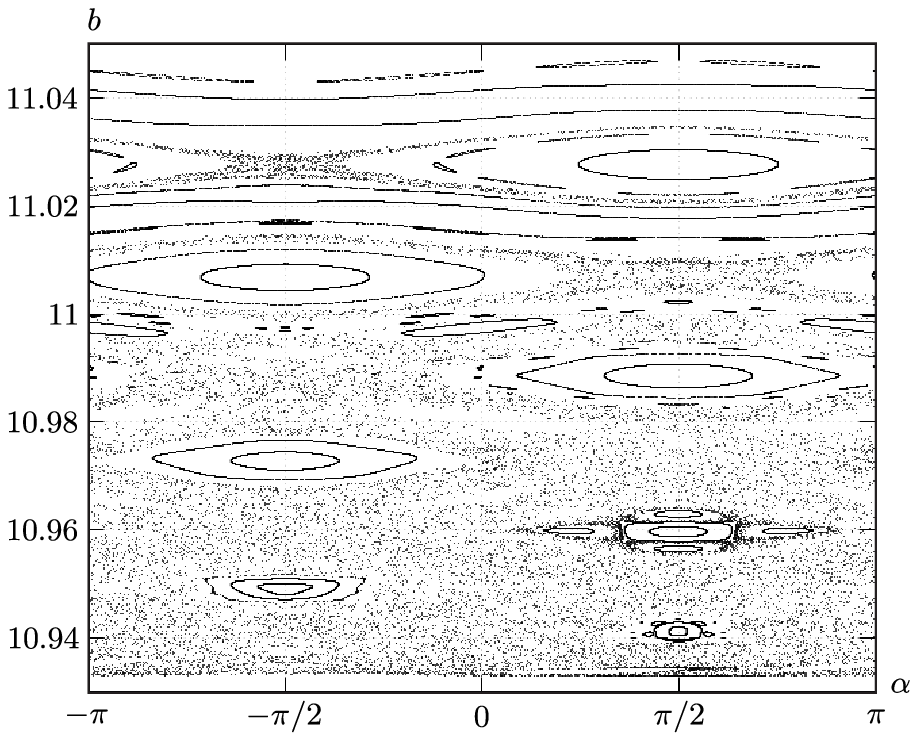}
	\caption{The parameter values $m_c = 1$, $\rho = 1$, $d = 0.01$, $I_c = 1$, $R = 1$, $k = 1$, $q = 1$, $h = 0.001$, $r_{max} = 100$. As the initial momentum we have chosen the largest value.} \label{fig.scatteringMap100max}
\end{figure}

\begin{figure}[h!]
	\centering
	\includegraphics{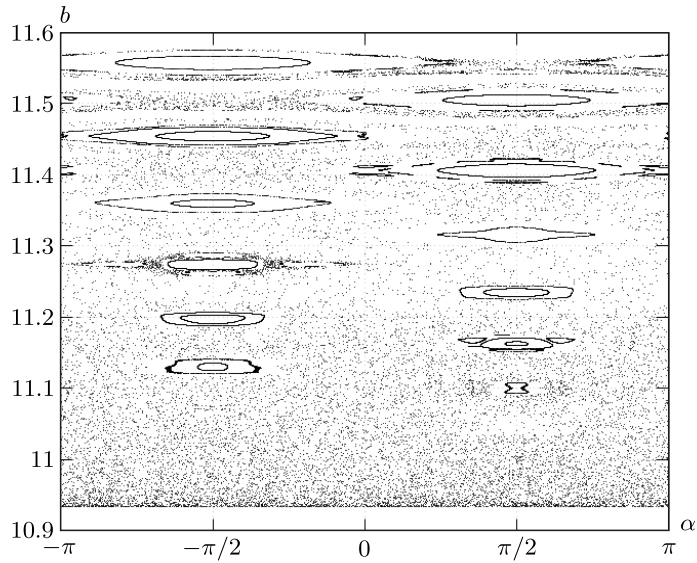}
	\caption{The parameter values $m_c = 1$, $\rho = 1$, $d = 0.01$, $I_c = 1$, $R = 1$, $k = 1$, $q = 1$, $h = 0.001$, $r_{max} = 100$. As the initial momentum we have chosen the smallest value.} \label{fig.scatteringMap100min}
\end{figure}

	We note that, by virtue of the properties \eqref{eq.asymptotic} of the reduced system,
the map \eqref{eq.scatteringMap} may be visualized more conveniently by passing from 
the angle $\varphi$ to the impact parameter:
	\begin{gather}
		b = r \sin(\alpha - \varphi) + m^{-1} m_c d \sin\alpha,
	\end{gather}
	whose value tends, as $r$ increases, to a constant value in contrast to the angle 
$\varphi$.
An example of a scattering map on the plane $(\alpha,\, b)$ is shown in Figs.~\ref{fig.scatteringMap100max} and \ref{fig.scatteringMap100min}

	Figures \ref{fig.scatteringMap100max} and \ref{fig.scatteringMap100min} shows invariant curves and chaotic regions similar to 
the classical Poincar\'{e} maps. In the case of a balanced circular foil ($d = 0$) the phase 
space of the scattering map is foliated by invariant curves $b = \const$.

\begin{figure}[h!]
		\centering
		\includegraphics{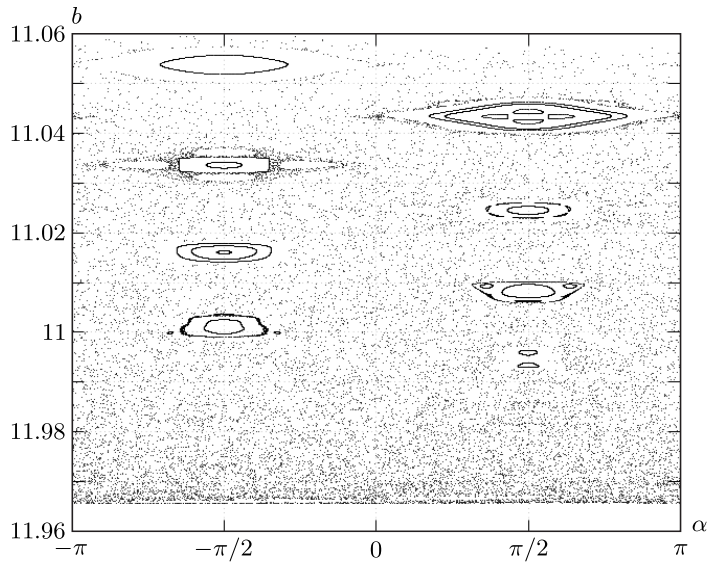}
		\caption{The parameter values $m_c = 1$, $\rho = 1$, $d = 0.01$, $I_c = 1$, $R = 1$, $k = 1$, $q = 1$, $h = 0.001$, $r_{max} = 200$. As the initial momentum we have chosen the largest value.} \label{fig.scatteringMap2}
	\end{figure}

	We note that, at large values of $r_{max}$, the portrait of the scattering map 
changes visually due to the drift of the variable $\alpha$ (see Fig. \ref{fig.scatteringMap2}).

	The question remains open whether nonintegrability with a given $r_{max}$
implies nonintegrability for any values of $r$. Thus, the constructed scattering maps 
allow the following hypothesis to be formulated:
	\begin{enumerate}
		\item[ ] \textit{The equations of motion \eqref{eq.equation} of an unbalanced
circular foil in the field of a fixed source of constant intensity are nonintegrable. }
	\end{enumerate}

	\begin{Remark}
		Constructing a scattering map requires a multiple integration of
equations \eqref{eq.subsys} by some numerical method for a long physical time interval. 
Since the system considered is Hamiltonian, the scattering map cannot contain any attracting 
or repelling sets. Nevertheless, when one uses explicit Runge\,--\,Kutta methods, 
the map may exhibit, in particular, asymptotically stable fixed points. This is due to 
the fact that long integration leads to a distortion of the level set of the energy integral 
\eqref{eq.Ham}. For this reason, for numerical integration we have used two types of 
methods: projection methods, which return the system's trajectory to the level surface
of the integral \eqref{eq.Ham} at the end of each step of integration, and collocation 
methods, which are symplectic and guarantee preservation of the phase volume with some 
accuracy \cite{Hairer_2006}. Both methods lead to visually similar results. However, 
projection methods are more efficient in terms of CPU time expenditures, since 
they can be based on any explicit method, whereas collocation methods are
completely implicit and require solving a nonlinear system of 
algebraic equations on each step of integration.
	\end{Remark}

	\section*{ACKNOWLEDGMENTS}
	
	The authors extend their gratitude to Prof. Oliver O'Reilly for useful comments. The authors are grateful to I. A. Bizyaev, I. S. Mamaev and A. A. Kilin for useful discussions.

	The work of Evgenii V. Vetchanin (introduction and section II) is supported by the RFBR under grant 18-29-10050-mk.
	The work of Elizaveta M. Artemova (sections I and III) was carried out within the framework of the state assignment of the Ministry of Education and Science of Russia (FEWS-2020-0009).

	\bibliography{biblio}		

\end{document}